\documentstyle[url,enumerate,graphicx,amsmath,a4]{isc6l}

\newcommand{\hX}{\widehat{X}}
\newcommand{\eqinlaw}{\stackrel{\mathcal{L}}{=}}

\newcommand{\E}[1]{\mathbf{E}\,#1}
\renewcommand{\Re}[1]{\mathrm{Re}\,#1}
\renewcommand{\Pr}[1]{\mathbf{P}\,#1}

\renewcommand{\mathbf}[1]{{\bf #1}}
\renewcommand{\mathcal}[1]{{\cal #1}}
\renewcommand{\mathrm}[1]{{\rm #1}}

\begin{document}

\begin{flushleft}
  {\LARGE\bf The space requirement of $m$-ary search trees:\ 
    distributional asymptotics for $m \geq 27$}

  \vspace{1.0cm}
  
  James Allen Fill$^1$\footnote{Research supported by \textsc{NSF}
    grant \textsc{DMS}--0104167 and by The Johns Hopkins University's
    Acheson J.~Duncan Fund for the Advancement of Research in
    Statistics.} and Nevin Kapur$^2$\footnote{Research partially
    supported by \textsc{NSF} grant~0049092.}

  \begin{description}
  \item $^1 \;$ Applied Mathematics and Statistics, The
    Johns Hopkins University, 3400~N.~Charles St., Baltimore MD
    21218-2682
  \item $^2 \;$ Computer Science, California Institute of
    Technology, MC~256-80, 1200~E.California Blvd., Pasadena CA 91125
  \end{description}
\end{flushleft}

\vspace{0.75cm}

\noindent {\bf Abstract}.
  We study the space requirement of $m$-ary search trees under the
  random permutation model when $m \geq 27$ is fixed.  Chauvin and
  Pouyanne have shown recently that $X_n$, the space requirement of an
  $m$-ary search tree on $n$~keys, equals $\mu (n+1) + 2
  \Re{[ \Lambda n^{\lambda_2}]} + \epsilon_n n^{\Re{\lambda_2}}$, where
  $\mu$ and $\lambda_2$ are certain constants, $\Lambda$ is a
  complex-valued random variable, and $\epsilon_n \to 0$ a.s.\ and
  in~$L^2$ as $n \to \infty$.  Using the contraction method, we
  identify the  distribution of~$\Lambda$.

\vskip 2mm


\section{Introduction}
\label{sec:introduction}

We start by giving a brief overview of search trees, which are
fundamental data structures in computer science used in searching and
sorting.  For integer~$m \geq 2$, the $m$-ary search tree, or multiway
tree, generalizes the binary search tree.  The quantity~$m$ is called
the \emph{branching factor}.  According to~\cite{MR93f:68045}, search
trees of branching factors higher than 2 were first suggested by Muntz
and Uzgalis~\cite{muntz71:_dynam} ``to solve internal memory problems
with large quantities of data.''  For more background we refer the
reader to~ \cite{knuth97,knuth98} and~\cite{MR93f:68045}.

An \emph{$m$-ary tree} is a rooted tree with at most $m$
``children'' for each \emph{node (vertex)}, each child of a node being
distinguished as one of $m$ possible types.  Recursively
expressed, an $m$-ary tree either is empty or consists of a
distinguished node (called the \emph{root}) together with an ordered
$m$-tuple of \emph{subtrees}, each of which is an $m$-ary
tree.

An \emph{$m$-ary search tree} is an
$m$-ary tree in which each node has the capacity to contain
$m-1$~elements of some linearly ordered set, called the set of
\emph{keys}.
In typical implementations of $m$-ary search trees, the keys at each
node are stored in increasing order and at each node one has $m$
pointers to the subtrees. By spreading the input data in $m$
directions instead of only 2, as is the case for a binary search tree,
one seeks to have shorter path lengths and thus quicker searches.

\begin{figure}[htbp]
  \centering
  \includegraphics[width=\textwidth]{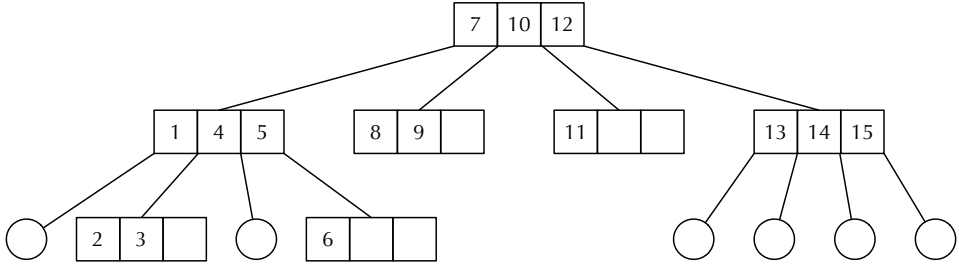}
  \caption{An $m$-ary search tree with space requirement~13.}
  \label{fig:m-ary}
\end{figure}

We consider the space of $m$-ary search trees on $n$ keys, and
assume that the keys are linearly ordered.  Hence, without loss
of generality, we can take the set of keys to be $[n] :=
\{1,2,\ldots,n\}$. We construct an $m$-ary search tree from a
sequence $s$ of $n$ distinct keys in the following way:
\begin{enumerate}[(i)]
\item If $n < m$, then all the keys are stored in the root node in
  increasing order.\label{item:4}
\item If $n \geq m$, then the first $m-1$ keys in the sequence
  are\label{item:5} stored in the root in increasing order, and the
  remaining $n-(m-1)$ keys are stored in the subtrees subject to the
  condition that if $\sigma_1 < \sigma_2 < \cdots < \sigma_{m-1}$
  denotes the ordered sequence of keys in the root, then the keys in
  the $j$th subtree are those that lie between $\sigma_{j-1}$ and
  $\sigma_{j}$, where $\sigma_0 := 0$ and $\sigma_{m} := n+1$,
  sequenced as in $s$.
\item All the subtrees are $m$-ary search trees that satisfy
  conditions~(\ref{item:4}),~(\ref{item:5}),
  and~(\ref{item:6}).\label{item:6}
\end{enumerate}
For example the $m$-ary search constructed from the sequence
\[
(10, 7,
12, 4, 1, 8, 5, 6, 9, 14, 11, 2, 15, 13, 3)
\]
is show in Figure~\ref{fig:m-ary}.  Note that empty nodes (also called
\emph{external nodes}) are represented as circles in the figure;
$m$~such nodes arise as children of a given node when that node
becomes filled to its capacity of $m-1$~keys.  In this paper
the total number of nodes (empty and nonempty) in an $m$-ary search
tree is called the \emph{space requirement} of the tree.

The uniform distribution on the space of permutations
of~$[n]$ induces a distribution of the space of $m$-ary search trees
with $n$~keys.  This is known as the \emph{random permutation model}.

Several authors have studied the limiting distribution of the space
requirement under the random permutation model.  Mahmoud and
Pittel~\cite{MR90a:68012} showed that when $m \leq 15$, the limiting
distribution is normal.  The result was later extended
to include $m \leq 26$ by Lew and Mahmoud~\cite{MR95i:68030}.  Chern
and Hwang~\cite{MR1871558} proved that when $m \geq 27$, the space
requirement centered by its mean and scaled by its standard deviation
does not have a limiting distribution.  Our result, stated as
Theorem~\ref{thm:d2}, for the case $m \geq 27$ was inspired by a recent
development (stated at the beginning of Section~\ref{sec:summary}) of
Chauvin and Pouyanne~\cite{Chauvin-Pouyanne-2003}.

\section{Summary}
\label{sec:summary}

Let $X_n$ denote the space requirement of an $m$-ary search tree on
$n$~keys chosen under the random permutation model.  Recently, Chauvin
and Pouyanne~\cite{Chauvin-Pouyanne-2003} have used martingale
techniques to show that when $m \geq 27$, we have $X_n = \hX_n
+ n^\sigma \epsilon_n$, where
\begin{equation}
  \label{eq:16}
  \hX_n := \frac{1}{H_m - 1} (n+1) + 2 \Re[n^{\lambda_2} \Lambda],
\end{equation}
with $\Lambda$ some complex-valued random variable and $\epsilon_n \to
0$ a.s.\ and in~$L^2$.  [In fact, they derive the asymptotics of the
random vector $(S_n^{(0)}, \dots, S_n^{(m-1)})$, where $S_n^{(i)}$
denotes the number of nodes with $i$ keys in a tree with $n$~keys, but
we shall be content here to study $X_n = \sum_{i=0}^{m-1} S_n^{(i)}$.]
In this representation, $\lambda_2 = \sigma + i \tau$ is the root of
the polynomial
\begin{equation}
  \label{eq:2}
  \phi(z) \equiv \phi_m(z) := (z+1) \cdots (z+m-1) - m!
\end{equation}
having second-largest real part and positive imaginary part.  It is
our goal to describe the distribution of the random variable~$\Lambda$.

To begin, we define the following
distributional transform~$T$ on $\mathcal{M}_2(\mu)$, the space of
probability distributions with a certain mean~$\mu$ defined
at~\eqref{eq:9} and finite second absolute moment:
\begin{equation}
  \label{eq:10}
  T:\ \mathcal{M}_2(\mu) \to \mathcal{M}_2(\mu), \quad
  \mathcal{L}(W) \mapsto \mathcal{L}\left( \sum_{k=1}^m
  S_k^{\lambda_2} W_k \right),
\end{equation}
where $(W_k)_{k=1}^m$ are independent copies of~$W$.  Here $\mathbf{S}
\equiv (S_1, \dots, S_m)$ is the vector of spacings of $m-1$
independent $\text{Uniform}(0,1)$ random variables~$U_1, \dots,
U_{m-1}$; i.e., if $U_{(1)}, \dots, U_{(m-1)}$ are their order
statistics and~$U_{(0)} := 0$, $U_{(m)} := 1$, then
\begin{equation}
  \label{eq:11}
  S_j := U_{(j)} - U_{(j-1)},\quad j=1,\dots, m.
\end{equation}
Furthermore, we take~$\mathbf{S}$ to be independent
of~$(W_k)_{k=1}^m$.
Next, define the metric~$d_2$ on~$\mathcal{M}_2(\mu)$ by
\begin{equation*}
  d_2(F,G) := \min \{\Vert X - Y \Vert_2:\ \mathcal{L}(X) = F,\,
  \mathcal{L}(Y) = G \},
\end{equation*}
with $\Vert X \Vert_2 := (\E{|X|^2})^{1/2}$ denoting the $L^2$-norm.
In the sequel, for notational convenience we will write $d_2(X, Y)$
instead of $d_2( \mathcal{L}(X), \mathcal{L}(Y) )$.

Our main result is the following.  (See the remark below
Lemma~\ref{lemma:a_J_k-b_J_k} for a strengthening.)
\begin{theorem}
  \label{thm:d2}
  Let $X_n$ denote the space requirement of an $m$-ary search tree on
  $n$~keys under the random permutation model with $m \geq 27$.  Define
  \begin{equation*}
    V_n := X_n - \frac{1}{H_m-1}(n+1)
  \end{equation*}
  and $\widehat{V}_n := 2 \Re{[n^{\lambda_2} Y]}$.  Here~$Y$ is a
  random variable with distribution equal to the unique fixed
  point~$\mathcal{L}(Y)$ of the distributional
  transform~\eqref{eq:10}.  Then
  $d_2(V_n, \widehat{V}_n) =  o(n^\sigma)$ and consequently
  $\Lambda$ has the same distribution  as $Y$.
\end{theorem}
The proof of Theorem~\ref{thm:d2} is presented in
Section~\ref{sec:proofs}, with the existence of the unique fixed point
established in Section~\ref{sec:exist-limit-distr} and bounds on
the $d_2$-distance derived in Section~\ref{sec:d_2-bounds}.
\begin{remark}
  \label{rem:urn-models}
  As discussed in~\cite{Chauvin-Pouyanne-2003} and~\cite{MR2040966},
  the study of the random vector~$(S_n^{(0)}, \dots, S_n^{(m-1)})$ can
  be recast as a generalized P{\'o}lya urn scheme which in turn can be
  studied by embedding into a continuous-time Markov multitype
  branching process.  Janson~\cite{MR2040966} obtains asymptotic
  distributional results for a very general class of urn schemes and
  multitype branching processes.  These include results for $m$-ary
  search trees, with~\eqref{eq:16} as a notable example.  We anticipate
  that our contraction-method technique for identifying
  $\mathcal{L}(\Lambda)$ in~\eqref{eq:16} will extend quite generally
  to oscillatory cases of Janson's results; this is the subject of
  ongoing research.  
  \hfill \qed
\end{remark}

In the sequel we will use $1 =: \lambda_1, \lambda_2, \dots,
\lambda_{m-1}$ to denote the $m-1$ roots of~\eqref{eq:2} in
nonincreasing order of real parts and roots with positive imaginary
parts listed before their conjugates.  In~\cite[\S3.3]{MR93f:68045}
and~\cite{FK-transfer}, the polynomial~$\psi(\lambda) =
\phi(\lambda-1)$ is considered.  The properties of the roots of~$\phi$
that we employ follow immediately from those known for the roots
of~$\psi$.

\section{Proofs}
\label{sec:proofs}

As preliminaries, note that the space requirement~$X_n$ has initial
conditions {$X_0 = X_1 = \cdots = X_{m-2} = 1$}, and for $n \geq
m-1$ that the number of keys \emph{not} stored in the root is
\begin{equation*}
  n' := n - (m-1).
\end{equation*}
It is well known that, under the random permutation model, $X_n$
satisfies the distributional recurrence 
\begin{equation}
  \label{eq:1}
  X_n \eqinlaw \sum_{k=1}^m X_{J_k}^{(k)} + 1, \quad n \geq m-1,
\end{equation}
where $\eqinlaw$ denotes equality in law (i.e.,\ in distribution), and
where, on the right, 
\begin{itemize}
\item the random vector~$\mathbf{J} \equiv (J_1, \dots, J_m)$ is
  uniformly distributed over all $m$-tuples $(j_1, \ldots, j_m)$ of
  nonnegative integers with $j_1 + \dots + j_m = n'$;
\item for each $k = 1, \dots, m$, we have $X_j^{(k)} \eqinlaw X_j$;
\item the quantities $\mathbf{J}; X_0^{(1)}, \dots, X_{n'}^{(1)};
  X_0^{(2)}, \dots, X_{n'}^{(2)}; \dots ; X_0^{(m)}, \dots,
  X_{n'}^{(m)}$ are all independent.
\end{itemize}

Using~\eqref{eq:1}, we get a distributional recurrence for~$V_n$, with
notation as for the $X$'s:
\begin{equation}
  \label{eq:12}
  V_n \eqinlaw\sum_{k=1}^m V_{J_k}^{(k)}, \quad n \geq m-1.
\end{equation}
The initial conditions here are $V_j = 1 - \frac{j+1}{H_m - 1}$ for
$j=0,1, \dots, m-2$.  The asymptotics of the mean of~$V_n$ can be
derived using~\cite[Equation~(2.7)]{FK-transfer}:
\begin{equation}
  \label{eq:9}
  \E{V_n} = \mu n^{\lambda_2} + \bar{\mu} n^{\lambda_3} +
  O(n^{\Re{\lambda_4}}),
\end{equation}
where $\mu$ is a constant.  Note that no two roots of~\eqref{eq:2}
have the same real part unless they are mutually conjugate, so that
$\Re{\lambda_4} < \Re{\lambda_3} = \Re{\lambda_2} = \sigma$.

For the reader's convenience, we state here a part of the Asymptotic
Transfer Theorem of~\cite{FK-transfer}.  We will use this result in
Section~\ref{sec:d_2-bounds}.  The constant~$K'$ can be expressed in
terms of~$K$, but we shall have no use here for such an expression.
\begin{proposition}
  \label{thm:att}
  For fixed $m \geq 2$, consider the recurrence
  \begin{equation*}
    a_n = b_n + \frac{m}{\binom{n}{m-1}} \sum_{j=0}^{n'}
    \binom{n-1-j}{m-2} a_j, \quad n \geq m-1,
  \end{equation*}
  with specified initial conditions $(a_j)_{j=0}^{m-2}$.  If $b_n =
  Kn^v + o(n^v)$ with $v > 1$ and $K$ a constant, then 
  \begin{equation*}
    a_n = K' n^v + o(n^v)
  \end{equation*}
  where $K'$ is a constant.
\end{proposition}

\subsection{Fixed point}
\label{sec:exist-limit-distr}
The existence and uniqueness of the fixed point of the map~$T$
at~\eqref{eq:10}
follows from the contraction method~(see, e.g.,~\cite{MR2003c:68096}).
Indeed a routine modification of the argument presented
in~\cite[\S6]{FK-transfer} yields that~$T$ is a contraction
on~$\mathcal{M}_2(\mu)$ with contraction factor
\begin{equation*}
  \rho = \left[ m! \frac{\Gamma(2\sigma + 1)}{\Gamma(2\sigma + m)}
  \right]^{1/2} = \left[ \frac{m!}{(2\sigma + m - 1) \cdots (2\sigma +
  1)} \right]^{1/2} < 1,
\end{equation*}
since for $m \geq 27$, we have~$\sigma >
1/2$~\cite{MR93f:68045,FK-transfer}.  

\subsection{$d_2$ bounds}
\label{sec:d_2-bounds}

We begin by defining $d_n := d_2(V_n, \widehat{V}_n)$ and $f(t) :=
2\,\Re t = t + \bar{t}$.  Unless otherwise noted we will henceforth
assume $n \geq m-1$.  Throughout $\sum_{\mathbf{j}}$ will denote a
sum over all $m$-tuples~$(j_1, \dots, j_m)$ of nonnegative integers
summing to~$n'$.

By the triangle inequality,
\begin{equation}
  \label{eq:6}
  d_n \leq a_n + b_n,
\end{equation}
where, taking $(Y_k)_{k=1}^m$ to be independent copies of the random
variable~$Y$ in Theorem~\ref{thm:d2} and $\mathbf{J}$ and~$\mathbf{S}$
each independent of $(Y_k)_{k=1}^m$,
\begin{equation}
  \label{eq:8}
  a_n := d_2\left( V_n, \sum_{k=1}^m f( J_k^{\lambda_2} Y_k) \right)
\end{equation}
and
\begin{equation}
  \label{eq:14}
  b_n := d_2\left( \sum_{k=1}^m f(J_k^{\lambda_2} Y_k), \sum_{k=1}^m f(
  n^{\lambda_2} S_k^{\lambda_2} Y_k) \right).
\end{equation}
We proceed by deriving upper bounds for $a_n$ and~$b_n$ separately.
The bound on~$b_n$ is proved as Lemma~\ref{lemma:b_n}.

For~$a_n$ a crude bound can be derived as follows.  Even though this
bound is not sufficient to show that $d_n = o(n^\sigma)$, it will be
employed in Lemma~\ref{lemma:b_J_k2}, which in turn will be used to
derive the estimate that we need.
\begin{lemma}
  \label{lem:a_n-trivial}
  With $a_n$ defined at~\eqref{eq:8},
  \begin{equation*}
    a_n = O(n^\sigma).
  \end{equation*}
\end{lemma}
\begin{proof}
  By the triangle inequality,
  \begin{equation*}
    a_n \leq \Vert V_n \Vert_2 + \sum_{k=1}^m \Vert
    f(J_k^{\lambda_2}Y_k) \Vert_2 = \Vert V_n \Vert_2 + m \Vert
    f(J_1^{\lambda_2} Y_1) \Vert_2.
  \end{equation*}
  Since $J_1 \leq n'$ and $\Vert Y_1 \Vert_2 < \infty$, we have
  $\Vert f(J_1^{\lambda_2} Y_1) \Vert_2 = O(n^\sigma)$.
  Using independence of the $V^{(k)}_{j_k}$'s, \eqref{eq:12},
  and~\eqref{eq:9}, we have
  \begin{align*}
    \Vert V_n \Vert_2^2 &= \sum_{\mathbf{j}} \Pr[ \mathbf{J} =
    \mathbf{j}]\, \E \left| \sum_{k=1}^m V_{j_k}^{(k)} \right|^2 =
    \frac{1}{\binom{n}{m-1}} \sum_{\mathbf{j}}
    \sum_{k=1}^m \Vert V_{j_k}
      \Vert_2^2 + O(n^{2\sigma}) \\
    &= \frac{m}{\binom{n}{m-1}} \sum_{j=0}^{n-(m-1)}
    \binom{n-1-j}{m-2} \Vert V_j \Vert_2^2 + O(n^{2\sigma}).
  \end{align*}
  It follows from Theorem~\ref{thm:att} that $\Vert V_n \Vert_2^2 =
  O(n^{2\sigma})$, and the result follows. \hfill \qed
\end{proof}

To sharpen Lemma~\ref{lem:a_n-trivial}, we employ the following
coupling between the distributions of~$V_n$ and of $\sum_{k=1}^m
f(J_k^{\lambda_2} Y_k)$.  The $L^2$~distance exhibited by this
coupling serves as an upper bound on the $d_2$-distance.  For $k=1,
\dots, m$, let $(V_1^{(k)}, V_2^{(k)}, \dots; Y_k)$ be independent
copies of $(V_1, V_2, \dots; Y)$ such that the coupling between $V_j$
and~$Y$ is $d_2$-optimal for each~$j$.  [To construct such a coupling,
first choose optimally-coupled $V_1$ and $Y$; having chosen $(V_1,
\dots, V_j; Y)$, choose $V_{j+1}$ so that it is optimally-coupled
with~$Y$.]  Then, with $\mathbf{J} \equiv (J_k)_{k=1}^m$ independent
of everything else,
\begin{equation}
   \label{eq:5}
  a_n^2 \leq \left\Vert \sum_{k=1}^m V_{J_k}^{(k)} - \sum_{k=1}^m
  f(J_k^{\lambda_2} Y_k) \right\Vert_2^2
  = \sum_{\mathbf{j}} \Pr[ \mathbf{J} = \mathbf{j} ] \left\Vert
  \sum_{k=1}^m V_{j_k}^{(k)} - \sum_{k=1}^m f(j_k^{\lambda_2} Y_k)
  \right\Vert_2^2.
\end{equation}
Now
\begin{align}
  & \left\Vert \sum_{k=1}^m V_{j_k}^{(k)} - \sum_{k=1}^m f(j_k^{\lambda_2}
    Y_k)  \right\Vert_2^2 \notag \\
  &= \sum_{k=1}^m \bigl\Vert V_{j_k}^{(k)} - f(j_k^{\lambda_2} Y_k)
  \bigr\Vert_2^2 + \E{\sum_{1 \leq k \ne l \leq m} \bigl[
    V_{j_k}^{(k)} - f(j_k^{\lambda_2} Y_k)\bigr] \overline{\bigl[
    V_{j_l}^{(l)} - 
  f(j_l^{\lambda_2} Y_l)\bigr]}} \notag \\
   &= \sum_{k=1}^m d_{j_k}^2 + \sum_{1 \leq k \ne l \leq m} \E{\bigl[
     V_{j_k}^{(k)} - f(j_k^{\lambda_2} Y_k)\bigr]} \E{\overline{\bigl[
    V_{j_l}^{(l)} -  f(j_l^{\lambda_2} Y_l)\bigr]}} \label{eq:3}
\end{align}
If we choose the mean~$\E{Y}$ to be $\mu$,
it follows from~\eqref{eq:9} that $\E{\bigl[ V_n - f(n^{\lambda_2} Y)
  \bigr]} = O(n^{\Re{\lambda_4}})$.  It follows then that the second
sum in~\eqref{eq:3} is~$O(n^{2 \Re{\lambda_4}}) = o(n^{2\sigma})$
  uniformly in~$\mathbf{j}$.
  Thus, from~\eqref{eq:5} and~\eqref{eq:3},
\begin{equation}
  \label{eq:4}
  a_n^2 \leq \E{ \sum_{k=1}^m d_{J_k}^2} + r_n,
\end{equation}
where $r_n = o(n^{2\sigma})$.

Next, we proceed to bound~$b_n$.
\begin{lemma}
  \label{lemma:b_n}
  With $b_n$ defined at~\eqref{eq:14},
  \begin{equation*}
    b_n = o(n^\sigma).
  \end{equation*}
\end{lemma}
\begin{proof}
  We take $Y_1, \dots, Y_m$ to be independent copies of~$Y$ and
  $(\mathbf{J}, \mathbf{S})$ independent of~$Y_1, \dots, Y_m$.  The
  conditional distribution of~$\mathbf{J}$ given $\mathbf{S} =
  \mathbf{s} \equiv (s_1, \dots, s_m)$ is taken to be
  $\text{Multinomial}(n', \mathbf{s})$.  Indeed this yields the
  distribution of the vector of sizes of the subtrees rooted at the
  root of a random $m$-ary search tree~\cite{MR2000e:68073}.  Then
\begin{align*}
  b_n &\leq \biggl\Vert \sum_{k=1}^m f(J_k^{\lambda_2} Y_k) - \sum_{k=1}^m
  f(n^{\lambda_2} S_k^{\lambda_2} Y_k) \biggr\Vert_2 && \\
  &\leq \sum_{k=1}^m \bigl\Vert f(J_k^{\lambda_2} Y_k) - f(n^{\lambda_2}
  S_k^{\lambda_2} Y_k) \bigr\Vert_2 && \\
  &\leq 2 \sum_{k=1}^m \Bigl\Vert\, [J_k^{\lambda_2} - (n S_k)^{\lambda_2}] Y_k
  \,\Bigr\Vert_2 && \text{(by definition of~$f$)}\\
  &= 2 \Vert Y \Vert_2 \, \sum_{k=1}^m \bigl\Vert J_k^{\lambda_2} - (n
  S_k)^{\lambda_2} \bigr\Vert_2 && \text{(by independence)} \\
  &= 2 m \Vert Y \Vert_2 \, \bigl\Vert J_1^{\lambda_2} - (n
  S_1)^{\lambda_2} \bigr\Vert_2. && \text{(by symmetry)}
\end{align*}
We know that $\Vert Y \Vert_2 < \infty$, and by
Lemma~\ref{lemma:j_k-nS_k} to follow the last factor above
is~$o(n^\sigma)$. \hfill\qed
\end{proof}

\begin{lemma}
  \label{lemma:j_k-nS_k}
  With $\sigma > 1/2$ denoting $\Re\lambda_2$,
  \begin{equation*}
    \Vert J_1^{\lambda_2} - (n S_1)^{\lambda_2} \Vert_2 =
    o(n^{\sigma}).
  \end{equation*}
\end{lemma}
\begin{proof}
  Given $\epsilon > 0$ we will show that the $L_2$-norm in question is
  bounded by a constant times $\epsilon^{1/2} n^\sigma$.   The lemma
  then follows by letting  $\epsilon \downarrow 0$.
  
  Observe that
  \begin{equation}
    \label{eq:13}
    \Vert J_1^{\lambda_2} - (n S_1)^{\lambda_2} \Vert_2^2 =
    \E{|J_1^{\lambda_2} - (n S_1)^{\lambda_2}|^2} = \E{
    \E\bigl[\,|J_1^{\lambda_2}
    - (n S_1)^{\lambda_2} |^2 \mid S_1\bigr]}.
  \end{equation}
  Until further notice assume $s > 2\epsilon$, and note that the
  conditional expectation $\E[\,| J_1^{\lambda_2} - (n
  S_1)^{\lambda_2} |^2 \mid S_1 = s]$ equals
  \begin{equation*}
    \sum_{j=0}^{n'} \Pr[ J_1 = j | S_1 = s] |j^{\lambda_2} - (n
    s)^{\lambda_2}|^2 
     = \sum_{0 \leq j \leq n(s - \epsilon)} + \sum_{
      n(s - \epsilon) < j < n(s+\epsilon)} +
    \sum_{n(s+\epsilon) \leq j \leq n}.
  \end{equation*}
  
  The conditional distribution of~$J_1$ given~$S_1 = s$ is
  $\text{Binomial}(n',s)$.  The last sum on the right is~$o(1)$
  uniformly in~$s$ since, by~\cite[Ex.~1.2.10-21]{knuth97},
  \begin{equation*}
    \Pr[J_1 \geq n(s + \epsilon)
    \mid S_1 = s] \leq \Pr[J_1 \geq n'(s + \epsilon)
    \mid S_1 = s] \leq \exp{(-\epsilon^2n'/2)}.
  \end{equation*}
  For the first sum observe that, for $n$ large enough (independently
  of~$s$),
  \begin{equation*}
    \Pr[J_1 \leq n(s-\epsilon) \mid S_1 = s] \leq \Pr\left[ J_1 \leq
    n'\left(s-\frac\epsilon2\right) \Bigl\vert S_1 = s\right] \leq
    \exp{(-\epsilon^2n'/8)},
  \end{equation*}
  the last inequality being a consequence of the aforementioned
  exercise.  Thus the first sum is also~$o(1)$ uniformly in~$s$.
  
  On the other hand, for the range of summation in the middle sum, by
  the mean value theorem and the assumed inequality $\epsilon < s/2$
  we have
  \begin{equation*}
    \left| \left( \frac{j}{n} \right)^{\lambda_2} - s^{\lambda_2} \right| \leq
      \epsilon |\lambda_2| \max_{\zeta \in (s-\epsilon, s+\epsilon)}
      |\zeta|^{\sigma-1} \leq \epsilon |\lambda_2| c_\sigma s^{\sigma-1},
  \end{equation*}
  where $c_\sigma$ is $(3/2)^{\sigma-1}$ if $\sigma \geq 1$ and
  $(1/2)^{\sigma-1}$ if $\sigma < 1$.  Thus
  \begin{equation*}
    |j^{\lambda_2} - (ns)^{\lambda_2}|^2 = n^{2 \sigma}\left| \left(
    \frac{j}{n} 
    \right)^{\lambda_2} - s^{\lambda_2} \right|^2 \leq \epsilon^2
    |\lambda_2|^2  c_\sigma^2 s^{2(\sigma-1)} n^{2\sigma}.
  \end{equation*}
  Hence the middle sum is at most $\epsilon^2 |\lambda_2|^2 c_\sigma^2
  s^{2(\sigma-1)} n^{2\sigma}$.

  Note that $S_1$ has distribution
  Beta$(1,m)$ and that
  \begin{equation*}
    \int_0^1 s^{2(\sigma-1)} (1-s)^{m-1}\,ds = \frac{\Gamma(m)
    \Gamma(2\sigma -1)}{\Gamma(m + 2\sigma - 1)} < \infty
  \end{equation*}
  since $\sigma > 1/2$.  So
  \begin{equation*}
    \int_{2\epsilon}^1 \E[\,| J_1^{\lambda_2} - (n S_1)^{\lambda_2}
    |^2 \mid S_1 = s]\, \Pr[ S_1 \in ds ] \leq \text{constant} \times
    \epsilon^2 n^{2\sigma}.
  \end{equation*}
  Finally,
  \begin{multline*}
    \int_0^{2\epsilon} \E[\,| J_1^{\lambda_2} - (n S_1)^{\lambda_2}
    |^2 \mid S_1 = s]\, \Pr[ S_1 \in ds ] \\
    \leq \text{constant} \times
    n^{2\sigma} \Pr[ S_1 \leq 2\epsilon] \leq \text{constant} \times
    \epsilon n^{2\sigma}.
  \end{multline*} \hfill\qed
\end{proof}

Combining~\eqref{eq:6} and~\eqref{eq:4}, we get
\begin{equation}
  \label{eq:7}
  a_n^2 \leq \E{ \sum_{k=1}^m (a_{J_k} + b_{J_k})^2} + r_n = \E{
  \sum_{k=1}^m a_{J_k}^2} + 2 \E{\sum_{k=1}^m a_{J_k} b_{J_k} } + \E{
  \sum_{k=1}^m b_{J_k}^2} + r_n.
\end{equation}
Next we bound the terms on the right-hand side, so that~\eqref{eq:7}
will yield a recursive inequality.
\begin{lemma}
  \label{lemma:b_J_k2}
  \begin{equation*}
    \E{ \sum_{k=1}^m b_{J_k}^2} = o(n^{2\sigma}).
  \end{equation*}
\end{lemma}
\begin{proof}
  By linearity of expectation and symmetry,
  \begin{equation*}
    \E{ \sum_{k=1}^m b_{J_k}^2} = \sum_{k=1}^m \E{b_{J_k}^2} = m\,
    \E{b_{J_1}^2}.
  \end{equation*}
  Now, the conditional distribution of~$J_1$ given~$S_1 =s$ is
  \mbox{$\text{Binomial}(n',s)$}.  We show that the
  conditional expectation $\E{[b_{J_1}^2 \mid S_1=s]}$ is
  $o(n^{2\sigma})$.  To that end,
  let $X$ be distributed $\text{Binomial}(n,s)$. For $\epsilon > 0$,
  \begin{equation*}
    \E{b_X^2} = 
    \sum_{j=0}^{n} \Pr[X = j ] b_j^2 = \sum_{0 \leq j
    \leq n(s-\epsilon)} +
    \sum_{ n(s - \epsilon) \leq j \leq n}.
  \end{equation*}
  Now an argument similar to the one used in the proof of
  Lemma~\ref{lemma:j_k-nS_k} can be employed.  The first sum on the
  right is~$o(n^{2\sigma})$.  On the other hand, we use the fact that
  $b_n = o(n^\sigma)$ from Lemma~\ref{lemma:b_n} to conclude that the
  second sum is $o(n^{2\sigma})$. \hfill\qed
\end{proof}
\begin{lemma}
  \label{lemma:a_J_k-b_J_k}
  \begin{equation*}
  \E{\sum_{k=1}^m a_{J_k} b_{J_k} } = o(n^{2\sigma}).
  \end{equation*}
\end{lemma}
\begin{proof}
  The proof (using the crude bound on~$a_n$ established in
  Lemma~\ref{lem:a_n-trivial}) is very similar to that of
  Lemma~\ref{lemma:b_J_k2}.  We omit the details. \hfill\qed
\end{proof}

We now complete the proof of Theorem~\ref{thm:d2}.
Using~\eqref{eq:7} and Lemmas~\ref{lemma:a_J_k-b_J_k}
and~\ref{lemma:b_J_k2} we find
\begin{align*}
  a_n^2 &\leq \E{\sum_{k=1}^m a_{J_k}^2} + g_n
  = \frac{1}{\binom{n}{m-1}} \sum_{\mathbf{j}} \sum_{k=1}^m a_{j_k}^2
  + g_n
  = \frac{m}{\binom{n}{m-1}} \sum_{\mathbf{j}} a_{j_1}^2 + g_n \\
  &= \frac{m}{\binom{n}{m-1}} \sum_{j=0}^{n-(m-1)} \binom{n-1-j}{m-2}
  a_j^2 + g_n,
\end{align*}
where $g_n = o(n^{2\sigma})$.  It follows from
Proposition~\ref{thm:att} that $a_n^2 = o(n^{2\sigma})$, so that
{$d_n \leq a_n + b_n = o(n^\sigma)$}, as desired.
\begin{remark}
  \label{rem:rates}
  The $o$-estimates in Lemmas~\ref{lemma:b_n}--\ref{lemma:a_J_k-b_J_k}
  can be improved to $O$-estimates.  In the proof of
  Lemma~\ref{lemma:j_k-nS_k}, choosing $\epsilon$ as a function of~$n$
  (specifically, taking $\epsilon_n$ to be a suitable constant
  multiple of $n^{-1/2} \log{n}$) sharpens the estimate $o(n^\sigma)$
  to~$O(n^{\sigma-\frac14}\sqrt{\log n})$, so that $b_n = O(n^{\sigma
    - \frac14} \sqrt{\log{n}})$ in Lemma~\ref{lemma:b_n}.  In turn,
  Lemmas~\ref{lemma:b_J_k2} and~\ref{lemma:a_J_k-b_J_k} are then
  immediately strengthened to $O(n^{2\sigma - \frac12} \ln{n})$ and
  $O(n^{2\sigma - \frac14} \sqrt{\log n})$, respectively.  This leads
  to $d_2(V_n, \widehat{V}_n) = O(n^{\Re{\lambda_4}}) + O(n^{\sigma -
    \frac18}(\log n)^{\frac14})$.  Numerics strongly support the
  conjecture that $\sigma - \Re{\lambda_4} \downarrow 0$ as $m
  \uparrow \infty$.  If this is true, then $d_2(V_n, \widehat{V}_n)$
  is $O(n^{\Re{\lambda_4}})$ whenever $m \geq 1044$.  Due to the
  presence of $r_n = O(n^{2 \Re{\lambda_4}})$ in~\eqref{eq:4}, this
  large-$m$~rate of convergence cannot be improved by the methods of
  this paper and presumably is the exact rate. \hfill\qed
\end{remark}

Finally, to prove equality in distribution of~$\Lambda$ and~$Y$, we
show that $d_2(\Lambda, Y) = 0$.  Indeed with $\Lambda =
|\Lambda|e^{i\Theta}$ and $Y = |Y| e^{iT}$, we have
\begin{align*}
  d_2\left( \Re(n^{\lambda_2} \Lambda), \Re(n^{\lambda_2}Y) \right) &=
  d_2 \left(\Re (n^{\sigma + i\tau} |\Lambda| e^{i\Theta} ), \Re
    (n^{\sigma + i\tau} |Y| e^{i T})\right) \\
  &= d_2\left( n^\sigma |\Lambda| \cos (\tau \ln n + \Theta), n^\sigma
  |Y| \cos(\tau \ln n + T) \right).
\end{align*}
But $d_2\left( \Re(n^{\lambda_2} \Lambda), \Re(n^{\lambda_2}Y)
\right) = o(n^\sigma)$ so that, as $n \to \infty$,
\begin{equation*}
  d_2 \left( |\Lambda| \cos(\tau \ln n + \Theta), |Y| \cos(\tau \ln n
  + T) \right) \to 0.
\end{equation*}
For any fixed~$\phi \in [0,2\pi)$ we can choose $n \to \infty$ such
that $(\tau \ln n) \bmod{2\pi} \to \phi$.  Then $|\Lambda|\cos(\phi +
\Theta)$ and $|Y|\cos(\phi + T)$ have the same distribution.  It
follows from the Cramer--Wold device~\cite[Theorem~29.4]{MR95k:60001}
that the random vectors $(|\Lambda| \cos\Theta, |\Lambda| \sin\Theta)$
and $(|Y| \cos{T}, |Y| \sin{T})$ have the same distribution.  In
particular, $\Lambda = |\Lambda| e^{i\Theta}$ and $Y = |Y| e^{iT}$
have the same distribution, as claimed.  This completes the proof
of Theorem~\ref{thm:d2}.

\bibliographystyle{habbrv}
\bibliography{msn,leftovers}

\end{document}